\documentclass{ws-p8-50x6-00}
\usepackage{amssymb}
\begin{document}
\newtheorem{theorem}{Theorem}[section]
\newtheorem{lemma}{Lemma}[section]
\newtheorem{remark}{Remark}[section]
\def\operatorname#1{{\rm#1\,}}
\def\text#1{{\hbox{#1}}}
\def\qedbox{\hbox{$\rlap{$\sqcap$}\sqcup$}}
\def\range{\operatorname{range}}
\def\Pspan{\operatorname{span}}
\def\rank{\operatorname{rank}}
\def\id{\operatorname{Id}}
\def\trace{\operatorname{trace}}
\font\pbglie=eufm10
\def\RR{\text{\pbglie R}}
\def\XJ{\text{\pbglie J}}
\def\AA{{\mathcal{A}}}
\def\CC{{\mathcal{C}}}
\def\JJ{{\mathcal{J}}}
\def\OO{{\mathcal{O}}}
\font\pbglie=eufm10
\def\SS{{\mathcal{S}}}
\def\pp{\phantom{-}}
\def\Gr{\text{Gr}}

 \makeatletter
  \renewcommand{\theequation}{%
   \thesection.\arabic{equation}}
  \@addtoreset{equation}{section}
 \makeatother
\title{Spacelike Jordan Osserman algebraic curvature tensors in the higher signature setting} 
\author{Peter B. Gilkey
 and Raina
Ivanova}
\address{Mathematics Department, University of Oregon, Eugene Or 97403 USA\\
email: gilkey@darkwing.uoregon.edu, ivanovar@hopf.uoregon.edu}
 \maketitle
\centerline{Dedicated to Professor Navaira}\bigbreak
\abstracts{Abstract: Let $R$ be an algebraic curvature tensor on a vector space of signature $(p,q)$
defining a spacelike Jordan Osserman Jacobi operator $\JJ_R$. We show
that the eigenvalues of
$\JJ_R$ are real and that $\JJ_R$ is diagonalizable if $p<q$.}

\section{Introduction}\label{sect1}
Let $V$ be a vector space that is equipped with an inner
product $(\cdot,\cdot)$ of signature
$(p,q)$. The inner product is said to be {\it
Riemannian} if
$p=0$, {\it Lorentzian} if
$p=1$, and {\it balanced} (or {\it neutral}) if
$p=q$. Let
$R$ be an algebraic curvature tensor on
$V$; i.e.
$R\in\otimes^4V^*$ satisfies the curvature symmetries of the Riemann curvature tensor:
\begin{eqnarray*}
&&R(x,y,z,w)=R(z,w,x,y)=-R(y,x,z,w),\text{ and}\\
&&R(x,y,z,w)+R(y,z,x,w)+R(z,x,y,w)=0.
\end{eqnarray*}
The associated {\it Jacobi operator}
$\JJ_R$ is the self-adjoint linear map of $V$ which is defined by the identity:
$$(\JJ_R(x)y,z)=R(y,x,x,z).$$
We say that an algebraic curvature tensor $R$ is {\it spacelike Osserman} if the eigenvalues of $\JJ_R$
are constant on the pseudo-sphere $S^+(V)$ of unit spacelike vectors. We say that $R$ is {\it spacelike
Jordan Osserman} if the Jordan normal form of
$\JJ_R$ is constant on $S^+(V)$. In the Riemannian context, these two notions are equivalent. However, in
higher signature setting, the eigenvalue structure does not determine the conjugacy class, so we shall
work with Jordan normal form rather than with the eigenvalues alone.

The investigation of spacelike Jordan Osserman tensors is motivated by geometric considerations. Let
$(M,g)$ be a pseudo-Riemannian manifold of signature
$(p,q)$. If
${}^gR$ is the Riemann curvature tensor of the Levi-Civita connection, then ${}^gR_P$ is an algebraic
curvature tensor on the tangent space
$T_PM$ for every $P$ in $M$. On the other hand, every algebraic curvature tensor is geometrically
realizable.
We say that $(M,g)$ is {\it spacelike Jordan Osserman} if the Jordan normal form
of $\JJ_{{}^gR}$ is constant on the pseudo-sphere bundle $S^+(M,g)$ of unit spacelike tangent vectors.

Let $(M,g)$ be Riemannian. Then 
$(M,g)$ is said to be a {\it local $2$ point homogeneous space} if the local isometries of $(M,g)$ act
transitively on $S^+(M,g)$; such a manifold is either flat or locally isometric to a rank $1$
symmetric space. If $(M,g)$ is a local $2$ point homogeneous space, then the Jacobi operator
$\JJ_R$ has constant Jordan normal form on
$S^+(M,g)$ and hence $(M,g)$ is spacelike Jordan Osserman.
Osserman\cite{refO} wondered if the converse held. Chi\cite{refChi} showed this was true if $m=4$, if
$m\equiv1$ mod $2$, or if $m\equiv2$ mod $4$; an important ingredient in Chi's work was the
classification of the spacelike Jordan Osserman Riemannian algebraic curvature tensors in these
dimensions. There are other partial results available\cite{refDD,refGSV}, but the question is still open
in the Riemannian context. 

It has been
shown\cite{refBBG,refGKV} that any Lorentzian spacelike
Jordan Osserman algebraic curvature tensor has constant sectional curvature. Thus
classification is complete in this context. 
Finally, it is known\cite{refBBGZ,refBCGHV,refGVV,refGrbc} that there exist balanced spacelike
Jordan Osserman pseudo-Riemannian manifolds which are not even locally homogeneous.

We shall work in the algebraic category henceforth. If $W$ is an auxiliary vector space and if $A$
is a linear map of $V$, then we define the {\it stabilization}
$$A\oplus 0:=\left(\begin{array}{ll}A&0\\0&0\end{array}\right)\text{ on }V\oplus W.$$

The Jordan normal form of a spacelike Jordan
Osserman algebraic curvature can be arbitrarily complicated in the balanced setting\cite{refGil,refGiIv}:

\begin{theorem}\label{arefa}
Let $J$ be a $r\times r$ real matrix and let $q\equiv0$ mod $2^r$. If $V$ is a vector space of neutral
signature
$(q,q)$, then there exists an algebraic curvature tensor $R$ on $V$ so that $\JJ_R(x)$ is conjugate to
$J\oplus0$ for every $x\in S^+(V)$.
\end{theorem}

However, the situation is very different if
$p<q$, i.e. if the spacelike directions in a certain sense dominate the timelike directions. We consider
this case and show that the geometry defined by such a tensor is much more rigid. The main result of this
paper is:

\begin{theorem}\label{arefb} Let $R$ be a spacelike Jordan Osserman algebraic curvature tensor on a vector
space $V$ of signature $(p,q)$, where $p<q$. Then $\JJ_R(x)$ is diagonalizable for any $x\in S^+(V)$.
\end{theorem}

Here is a brief outline of the paper. In \S\ref{sect2}, we review certain results concerning
self-adjoint maps in the indefinite setting. In \S\ref{sect3}, we establish several technical results
concerning vector bundles over projective spaces. In \S\ref{sect4}, we use the results of the previous
sections to prove Theorem
\ref{arefb}.

\section{Results from Linear Algebra}\label{sect2}

Let $\Re(\lambda)$ and $\Im(\lambda)$ be the real and imaginary parts of a
complex number $\lambda$. If $J$ is a linear map of a  real vector space $V$ of dimension $m$, then let
$J_\lambda$ be the real operator on $V$ defined by:
\begin{equation}J_\lambda:=\left\{\begin{array}{ll}
J-\lambda\cdot\id&\quad\text{if }\quad\lambda\in\mathbb{R},\\
(J-\lambda\cdot\id)(J-\bar\lambda\cdot\id)&\quad\text{if
}\quad\lambda\in\mathbb{C}-\mathbb{R}.\end{array}\right.\label{brefax}\end{equation}
We define the generalized eigenspaces by setting 
\begin{equation}E_\lambda=E_\lambda^J:=\ker\{J_\lambda^m\}.\label{brefay}\end{equation}

\begin{lemma}\label{brefa} Let $V$ be a vector space of signature $(p,q)$ and let $J$ be a self-adjoint
linear map of $V$. Then  $V$ can be decomposed as an orthogonal direct sum
$V=\oplus_{\Im(\lambda)\ge0}E_\lambda$. Furthermore, the induced metrics on the generalized eigenspaces
$E_\lambda$ are non-degenerate.
\end{lemma}

\medbreak\noindent{\bf Proof:} Let $\lambda$ and $\mu$ be complex numbers with $\lambda\ne\mu$
and
$\lambda\ne\bar\mu$. Since
$J_\lambda^m$ is self-adjoint and vanishes on
$E_\lambda$, we have
$$0=(J_\lambda^mx_\lambda,x_\mu)=(x_\lambda,J_\lambda^mx_\mu)\text{ for }x_\lambda\in E_\lambda\text{ and
} x_\mu\in E_\mu.$$
Since $J$ commutes with $J_\mu$, $J$ preserves $E_\mu$. Since the eigenvalues of
$J$ on $E_\mu$ are $\mu$ and $\bar\mu$, the linear maps $J-\lambda\cdot\id$, $J-\bar\lambda\cdot\id$, and
hence
$J_\lambda$ are isomorphisms of $E_\mu$; thus $J_\lambda^m(E_\mu)=E_\mu$. It now follows that
\begin{equation}E_\lambda\perp E_\mu\text{ and }
E_\lambda\cap E_\mu=\{0\}.\label{bork1}\end{equation}

 Let $V^{\mathbb{C}}:=V\otimes\mathbb{C}$ be the complexification of
$V$. We extend $J$ to $V^{\mathbb{C}}$ to be complex linear and set
$E_\lambda^{\mathbb{C}}:=\ker\{(J-\lambda)^m\}$. A complex vector space may be decomposed as the direct
sum of the generalized complex eigenspaces defined by a linear transformation. Consequently,
\begin{equation}V^{\mathbb{C}}=\oplus_\lambda E_\lambda^{\mathbb{C}}.\label{bork2}\end{equation}
As $E_\lambda^{\mathbb{C}}\oplus E_{\bar\lambda}^{\mathbb{C}}=E_\lambda\otimes\mathbb{C}$,
$V=\oplus_{\Im(\lambda)\ge0}E_\lambda$. By display (\ref{bork1}), the direct sum given in equation
(\ref{bork2}) is orthogonal; thus, the induced metric on each $E_\lambda$ is non-degenerate.
\qedbox

\medbreak We recall a few notions from bundle theory which we will use in what follows. Let
$\rho:E\rightarrow M$ be a vector bundle over a smooth manifold $M$. {\it The  fibers} $E_P:=\rho^{-1}(P)$
are real vector spaces which vary smoothly with the point $P\in M$. A {\it non-degenerate fiber metric} on
$E$ is a collection of non-degenerate inner products on each fiber which vary smoothly on $M$. A {\it bundle morphism} $\psi$ of $E$
is a collection of smooth linear maps $\psi_P$ of the fibers $E_P$ which vary smoothly with $P$. We say
$\psi$ is invertible if each $\psi_P$ is invertible. We say $\psi$ is {\it self-adjoint} if each $\psi_P$
is self-adjoint.

Let $V$ be a
vector space with a non-degenerate inner product. We can decompose
$V$ as a direct sum $V^+\oplus V^-$ of complementary orthogonal subspaces, where $V^+$ is
a maximal spacelike subspace and
$V^-$ is a maximal timelike subspace. There is a similar decomposition possible in the vector bundle
setting:

\begin{lemma}\label{brefB} Let $E$ be a vector bundle over a smooth manifold $M$ which is
equipped with a non-degenerate fiber metric. Then we can decompose $E$ as a direct sum
$E^+\oplus E^-$ of complementary orthogonal subbundles, where $E^+$ is a maximal spacelike subbundle and
$E^-$ is a maximal timelike subbundle.\end{lemma}

\medbreak\noindent{\bf Proof:} As noted above, we can decompose each individual fiber as an orthogonal
direct sum of a maximal spacelike and a maximal timelike subspace. The main technical difficulty is to
make the decompositions vary smoothly with
$P\in M$. We can use a partition of unity to put a positive definite inner product
$(\cdot,\cdot)_e$ on $E$. Define a bundle morphism $\psi$ of $E$ by setting
$(v,w)=(\psi v,w)_e$. Since each linear map $\psi_P$ is self-adjoint with respect to the positive
definite inner product
$(\cdot,\cdot)_e$ on each fiber $E_P$, $\psi_P$ is diagonalizable and has only real eigenvalues. As the
original inner product $(\cdot,\cdot)$ is non-degenerate,
each $\psi_P$ is invertable. Let $E_\lambda(\psi_P)\subset E_P$ be the eigenspaces of $\psi_P$ on
$E_P$. We define:
$$E^-_P:=\oplus_{\lambda<0}E_\lambda(\psi_P)\text{ and }
  E^+_P:=\oplus_{\lambda>0}E_\lambda(\psi_P).$$
By Lemma \ref{brefa}, the subspaces $E^+_P$ and $E^-_P$ are orthogonal and complementary. Since $\psi_P$
is invertible, the fibers $E^+_P$ and $E^-_P$ have constant rank and define smooth subbundles of $E$.
\qedbox

\medbreak We continue our preparation for the proof of Theorem \ref{arefb} by
studying vector bundles equipped with non-degenerate fiber metrics and self-adjoint bundle
morphisms which have constant Jordan normal form:

\begin{lemma}\label{brefC} Let $E$ be a vector bundle over a smooth manifold $M$ which is
equipped with a non-degenerate fiber metric. Let $J$ be a self-adjoint bundle morphism of $E$ which has
constant Jordan normal form. Let $\lambda$ be an eigenvalue of $J$. If
$J_\lambda\ne0$ on $E_\lambda$, choose $i\ge1$ maximal so
$J_\lambda^i(E_\lambda)\ne0$. Then $J_\lambda^iE_\lambda$ is a totally
isotropic subbundle of $E$ of non-zero rank.
\end{lemma}

\medbreak\noindent{\bf Proof:} Assume that $E$ and $J$ satisfy the hypothesis of the Lemma. Set:
$$E_{\lambda,i}:=J_\lambda^i(E_\lambda).$$
Since $J$ has constant Jordan normal form, $E_{\lambda,i}$ is a smooth vector bundle over $M$. Fix a
point $P$ of $M$ and let $v_1$ and $v_2$ be vectors in the fiber $E_{\lambda,i}(P)$. There exist vectors
$w_1,w_2\in E_{\lambda,i}(P)$ so
$v_1=J_\lambda^iw_1$ and
$v_2=J_\lambda^iw_2$. Note that $2i\ge i+1$, that
$J_\lambda$ is self-adjoint, and that $J_\lambda^{2i}=0$ on $E_\lambda$. We demonstrate that
$E_{\lambda,i}$ is totally isotropic and thereby complete the proof by computing:
$$(v_1,v_2)=(J_\lambda^iw_1,J_\lambda^iw_2)=(J_\lambda^{2i}w_1,w_2)=0.\ \qedbox$$

\section{Bundles over projective space}\label{sect3}

Let $V$ be a vector space of signature $(p,q)$. We decompose $V=V^+\oplus V^-$ as an orthogonal direct
sum, where
$V^+$ is a maximal spacelike subspace of dimension $q$ and $V^-$ is the complementary maximal timelike
subspace of dimension $p$. Let
$\mathbb{P}(V^+)$ be the projective space of lines in $V^+$. We define trivial bundles over
$\mathbb{P}(V^+)$ by setting:
\begin{eqnarray}
  &&\mathbb{V}^+:=\mathbb{P}(V^+)\times V^+,\nonumber\\ 
  &&\mathbb{V}^-:=\mathbb{P}(V^+)\times V^-,\text{ and}\label{crefx1}\\
  &&\mathbb{V}:=\mathbb{V}^+\oplus\mathbb{V}^-=\mathbb{P}(V^+)\times V.\nonumber\end{eqnarray}
These bundles inherit natural inner products from the given inner product on $V$. Let
$\langle x\rangle:=x\cdot\mathbb{R}$ be the line thru an element $x\in S(V^+)$.  The {\it classifying line
bundle}
$\gamma$ and the {\it orthogonal complement} $\gamma^\perp$ over
$\mathbb{P}(V^+)$ are the subbundles of $\mathbb{V}^+$ defined by:
\begin{eqnarray}
&&\gamma:=\{(\langle x\rangle,y)\in\mathbb{P}(V^+)\times V^+:y\in\langle x\rangle\},\text{ and}\nonumber\\
&&\gamma^\perp:=\{(\langle x\rangle,y)\in\mathbb{P}(V^+)\times V^+:y\perp\langle
x\rangle\}.\label{crefx2}\end{eqnarray}
Note that $\gamma_{\langle x\rangle}=\langle x\rangle$, i.e. the fiber of $\gamma$ over an element
$\langle x\rangle$ of $\mathbb{P}(V^+)$ is just the line
$\langle x\rangle$ itself. For that reason $\gamma$ has also been called the {\it tautological line
bundle}. This bundle plays an important role in the classification of real line bundles and in our further
considerations. In the following Lemma, we compare non-trivial (i.e. positive rank) subbundles of
$\gamma^\perp$ and $\mathbb{V}^-$.

\begin{lemma}\label{crefA} Let $V=V^+\oplus V^-$ be a vector space of signature $(p,q)$, where $p<q$.
Let $\mathbb{V}^-$ and $\gamma^\perp$  be the bundles over $\mathbb{P}(V^+)$ defined in equations
{\rm(\ref{crefx1})} and {\rm(\ref{crefx2})}, respectively. Then no non-trivial subbundle of
$\mathbb{V}^-$ is isomorphic to a non-trivial subbundle of
$\gamma^\perp$.\end{lemma}

\medbreak\noindent{\bf Proof:} The Stiefel-Whitney classes\cite{refMilnor} are cohomological
invariants of vector bundles.
  Let $w_1:=w_1(\gamma)$ be the first
Stiefel-Whitney class of the classifying line bundle $\gamma$. The cohomology ring of the
projective space $\mathbb{P}(V^+)$ is the truncated polynomial ring\cite{refMilnor}:
$$H^*(\mathbb{P}(V^+);\mathbb{Z}_2)=\mathbb{Z}_2[w_1]/(w_1^q=0).$$

To prove Lemma \ref{crefA}, we suppose the contrary, i.e. that there exist positive rank subbundles $E_1$
of
$\gamma^\perp$ and $E_2$ of $\mathbb{V}^-$ so that
$E_1$ is isomorphic to $E_2$, and argue for a contradiction. Let $r=\text{rank}(E_1)$.
Since $E_1\subset\gamma^\perp$, we may use Lemma 4.4.4\cite{refGil} to see that $w_r(E_1)=w_1^r$.
Furthermore, as
$E_1$ is isomorphic to $E_2$, we may conclude
$w_r(E_2)=w_1^r$. As
$\mathbb{V}^-=E_2\oplus E_2^\perp$, we have the following factorization
$$1=w(E_2)w(E_2^\perp)=(1+...+w_1^r)w(E_2^\perp)\quad\text{ in}\quad H^*(\mathbb{P}(V^+);\mathbb{Z}_2).$$
Since $\rank E_2=p-r<q$, the
truncation $(w_1^q=0)$ in $H^*(\mathbb{P}(V^+);\mathbb{Z}_2)$ plays no role, so we have the factorization
$$1=(1+...+w_1^r)w(E_2^\perp)\quad\text{ in }\quad\mathbb{Z}_2[w_1],$$
which is impossible. \qedbox

We use Lemma \ref{crefA} to establish the following Lemma:

\begin{lemma}\label{crefB} Let $V=V^+\oplus V^-$ be a vector space of signature $(p,q)$, where $p<q$.
Let $\mathbb{V}^-$ and $\gamma^\perp$  be the bundles over $\mathbb{P}(V^+)$ defined in equations
{\rm(\ref{crefx1})} and {\rm(\ref{crefx2})}, respectively. There is no totally isotropic non-trivial 
subbundle of 
$\gamma^\perp\oplus\mathbb{V}^-$.
\end{lemma}

\medbreak\noindent{\bf Proof:} We suppose, to the contrary, that there exists a  totally
isotropic non-trivial subbundle $E$ of $\gamma^\perp\oplus\mathbb{V}^-$. Let $\pi^+$ be orthogonal
projection on
$\gamma^\perp$ and let $\pi^-$ be orthogonal projection on $\mathbb{V}^-$. Set
$$E^+:=\pi^+(E)\subset\gamma^\perp\quad\text{and}\quad
  E^-:=\pi^-(E)\subset\mathbb{V}^-.$$
Note that $\ker\pi^+=\mathbb{V}^-$ and $\ker\pi^-=\gamma^\perp$. As $E$ is totally isotropic, every
vector in $E$ is null. Thus
$$E\cap\mathbb{V}^-=E\cap\gamma^\perp=\{0\}.$$ Consequently,
$\pi^+$ and $\pi^-$ define isomorphisms between $E$ and $E^+$ and between $E$ and $E^-$, respectively.
Thus $E^+$, which is a non-trivial subbundle of $\gamma^+$, is isomorphic to $E^-$, which is a
non-trivial subbundle of
$\mathbb{V}^-$. This contradicts Lemma
\ref{crefA}.
\qedbox

\section{Proof of Theorem \ref{arefb}}\label{sect4}

Let $V=V^+\oplus V^-$ be a vector space of signature $(p,q)$, where $p<q$, and
where $V^+$ and $V^-$ are orthogonal maximal spacelike and timelike subspaces of $V$,
respectively. Let $R$ be a spacelike Jordan Osserman algebraic curvature tensor on  $V$. Let $x\in
S^+(V)$. Since
$\JJ_R(x)$ is self-adjoint and since $\JJ_R(x)x=0$, $\JJ_R(x)$ preserves the orthogonal complement
$x^\perp$; we let
$\tilde\JJ_R(x)$ denote the restriction of $\JJ_R(x)$ to $x^\perp$; this is often called the {\it reduced
Jacobi operator}. The Jacobi operator can be represented in the form
$$\JJ_R(x)=\left(\begin{array}{cc}0&0\\\ 0&\ \tilde\JJ_R(x)\end{array}\right)\quad\text{on}\quad
x\oplus x^\perp.$$
Thus to prove Theorem
\ref{arefb}, it suffices to show that $\tilde\JJ_R(x)$ is diagonalizable.

If $\lambda\in\mathbb{C}$, then let $\tilde J_\lambda$ and
$E_\lambda$ be defined by $\tilde\JJ_R$ using equations (\ref{brefax}) and (\ref{brefay}),
respectively. We may then use Lemma \ref{brefa} to decompose
$$\gamma^\perp\oplus\mathbb{V}^-=\oplus_{\Im(\lambda)\ge0}E_\lambda
\text{ over }\mathbb{P}(V^+),$$
where the induced metric on each eigenbundle $E_\lambda$ is non-degenerate.
By Lemma \ref{crefB}, $E_\lambda$ does not contain a totally
isotropic subbundle. Thus, by Lemma \ref{brefC}, $\tilde\JJ_\lambda=0$ on $E_\lambda$. Consequently
$\tilde\JJ_R$ is diagonalizable on $E_\lambda$ if $\lambda\in\mathbb{R}$. 

To complete the proof, we must show that all the eigenvalues are real. Suppose, to the contrary, that
there exists an eigenvalue
$\lambda$ of
$\tilde\JJ_R$ so that
$\Im(\lambda)\ne0$; we argue for a contradiction. By
Lemma \ref{brefB}, $E_\lambda=E_\lambda^+\oplus E_\lambda^-$ decomposes as the orthogonal
direct sum of maximal spacelike and timelike subbundles. We define a bundle map
$\mathcal{I}$ of
$E_\lambda$ by setting:
$$\mathcal{I}:=\textstyle\frac{\tilde\JJ_R-\Re(\lambda)\id}{\Im(\lambda)}.$$
The definition of $\tilde\JJ_\lambda$ given in (\ref{brefax}) and the fact that
$\tilde\JJ_\lambda=0$ on
$E_\lambda$ then imply that
$\mathcal{I}^2=-\text{id}$ on $E_\lambda$.
Since $\mathcal{I}$ is self-adjoint, $\mathcal{I}$ is a para-isometry of $E_\lambda$ that
interchanges the roles of spacelike and timelike vectors. Thus,
$\mathcal{I}$ defines an isomorphism between $E_\lambda^+$ and $E_\lambda^-$. 

Let $\pi^+$ and $\pi^-$ be
orthogonal projections on $\gamma^\perp$ and $\Bbb{V}^-$, respectively. Since $E^+_\lambda$ contains no
timelike vectors and since $\ker(\pi^+)=\Bbb{V}^-$ is timelike, $\ker\pi^+\cap E_\lambda^+=\{0\}$ and
$\pi^+$ is an isomorphism from
$E^+_\lambda$ to
$\pi^+(E_\lambda)$. Similarly,
$\pi^-$ is an isomorphism from $E^-_\lambda$ to $\pi^-(E_\lambda)$. Thus $\pi^+(E_\lambda)$, which is
a non-trivial subbundle of $\gamma^\perp$, is isomorphic to $\pi^-(E_\lambda^-)$, that is a non-trivial
subbundle of
$\Bbb{V}^-$. This contradicts Lemma
\ref{crefA}.
\qedbox

\section*{Acknowledgments} Research of both authors supported in part by the NSF (USA) and MPI
(Leipzig).


\begin{thebibliography}{AAA}


\bibitem{refBBG} N. Bla\v zi\'c, N. Bokan and P. Gilkey,
    {\it A Note on Osserman Lorentzian manifolds}, Bull. London Math. Soc.
    {\bf 29} (1997), 227--230.

\bibitem{refBBGZ} N. Bla\v zi\'c, N. Bokan, P. Gilkey and Z. Raki\'c,{\it
     Pseudo-Riemannian Osserman manifolds},  J. Balkan Soc.
     Geometers {\bf l2} (1997), 1--12.

\bibitem{refBCGHV}
A. Bonome, R. Castro, E. Garc\'{\i}a--R\'{\i}o, L. Hervella, R. V\'{a}zquez--Lorenzo,
{\it Nonsymmetric Osserman indefinite K\"{a}hler manifolds},
 Proc. Amer. Math. Soc. {\bf 126} (1998), 2763--2769.


\bibitem{refChi} Q.-S. Chi, {\it A curvature characterization of certain locally
rank-one symmetric spaces},  J. Differential Geom. {\bf 28} (1988), 
     187--202.

\bibitem{refDD}I. Dotti and M. Druetta,
{\it Negatively curved homogeneous
     Osserman spaces},  Differential Geom. Appl. {\bf 11} (1999), 163--178.

\bibitem{refGKV} E. Garc\'ia-Ri\'o, D. Kupeli and M. E.
 V\'azquez-Abal, {\it On a problem of Osserman in Lorentzian geometry},
Differential Geom. Appl. {\bf 7} (1997), 85--100.

\bibitem{refGrbc}E. Garc\'ia-Rio, D. N. Kupeli, and R. V\' azquez-Lorenzo
{\bf Osserman manifolds in semi-Riemannian geometry}, Lecture notes in Mathematics,
Springer Verlag, to appear.

\bibitem{refGVV} E. Garc\'ia-Ri\'o, M. E. V\' azquez-Abal and
     R. V\' azquez-Lorenzo, {\it Nonsymmetric Osserman pseudo-Riemannian manifolds},
    Proc. Amer. Math. Soc. {\bf 126} (1998), 2771--2778.



\bibitem{refGil} P. Gilkey, {\bf Natural Operators Defined by the Riemann Curvature Tensor},
World Scientific Press, ISBN 981-02-4752.

\bibitem{refGiIv} P. Gilkey and R. Ivanova, {\it The Jordan normal form of Osserman algebraic
curvature tensors}, preprint.

\bibitem{refGSV} P. Gilkey, A. Swann, and L. Vanhecke,
     {\it Isoparametric geodesic spheres and a conjecture of Osserman regarding the
     Jacobi operator},  Quart. J. Math. Oxford Ser. {\bf 46} (1995),
      299--320.

\bibitem{refMilnor} J. Mil.nor and J. Stasheff, {\bf Characteristic Classes}, Annals of Math. Studies,
Princeton University Press (1974).

\bibitem{refO} R. Osserman, {\it Curvature in the eighties},  Amer. Math.
    Monthly {\bf 97} (1990), 731--756.
\end{thebibliography}
\end{document}